\documentclass[12pt]{amsart}

\usepackage{amsrefs}

\usepackage{xcolor}

\newtheorem*{theorem}{Theorem}

\theoremstyle{definition}
\newtheorem*{definition}{Definition}

\begin{document}
\title{A Hopf algebra without modular pair in
involution}
\author{Sebastian Halbig}
\author{Ulrich Kr\"ahmer}
\address{TU Dresden\\
Institut f\"ur Geometrie\\
01062 Dresden}
\email{sebastian.halbig@tu-dresden.de}
\email{ulrich.kraehmer@tu-dresden.de}
 
\begin{abstract}
The aim of this short note is to communicate
an example of a finite-dimensional Hopf
algebra that does not admit a modular pair in
involution in the sense of Connes and
Moscovici. 
\end{abstract}
\maketitle

\section{Introduction}
The concept of a modular pair in involution
was introduced by Connes and Moscovici
\cite{cm} in order to define the Hopf-cyclic
cohomology of a Hopf algebra $H$ over a field
$k$. 
In the following we freely use 
standard notation from Hopf
algebra theory e.g.~as in \cite{m,r}. 
In particular,  
$H^\circ$ is the Hopf dual of $H$ and 
$\beta ^{-1} =
\beta \circ S$ is the convolution inverse of
a group-like $ \beta \in H^\circ$ (i.e.~a 
character $ \beta \colon H \rightarrow k$). 

\begin{definition}
Let $H$ be a Hopf algebra. A
pair $(l, \beta )$ of group-like elements 
$l \in H, \beta \in H^\circ$ 
is a \emph{modular pair in
involution} if $ \beta (l) = 1$ and  
$$
	S^2(h) = \beta (h_{(1)}) l h_{(2)} l^{-1}
\beta ^{-1} (h_{(3)}) 
$$
holds for all $h \in H$. 
\end{definition}

Hajac et al.~extended this 
notion to that of stable anti
Yetter-Drinfel'd modules over Hopf algebras,
see \cite{hkrs}.
It is also related to earlier work by
Kauffman and Radford
\cite{kr} who classified the
ribbon elements in Drinfel'd doubles of
finite-dimensional Hopf algebras. Among
their results they showed that if $\dim _k H$
is odd and $S^2$ has odd order, then there is
always a pair $(l, \beta )$ implementing
$S^2$. The question arises whether there are
also always pairs $(l,\beta )$ that additionally satisfy the
stability condition $\beta (l) = 1$. 
The aim of the present note is to
point out that this is not the case in
general:

\begin{theorem}
Let $p$ be a prime number, $s \in
\mathbb{Z}_p \setminus \!\{0\}$, 
$q \in k$ be a primitive $p$-th root of
unity, and $H$ be the Hopf algbera with
generators $g,x,y$ and defining algebra and
coalgbera relations
$$
	gx = q x g ,\quad
	g y = q^{-s} y g,\quad
	g^p=1,\quad
	x^p = y^p =0,\quad
	xy = q^{-s} yx,
$$
$$
	\Delta (g) = g \otimes g,\quad
	\Delta (x) = 1 \otimes x + x \otimes
g,\quad
	\Delta (y) = 1 \otimes y +
	y \otimes g^s.
$$
Its antipode is determined by
$$
	S(g)=g^{-1}, \quad 
	S(x)= -xg^{-1} \quad
	S(y)= -yg^{-s},
$$
and $H$ has a modular pair in
involution if and only if 
$s \in \{0,1,p-1\}$.
\end{theorem}

The Hopf algebra $H$ appears naturally in
several contexts. In particular, it is
referred to as the book Hopf algebra in
\cite{as}. 

\subsection*{Acknowledgements} 
We thank P.M.~Hajac for pointing us to the
question answered here. 

\section{Proof}
It is immediately verified that the
group-likes in $H$ are the elements of the form 
$l=g^i$ for some $i$; furthermore, 
a character $ \beta \colon H \rightarrow k$,
has to vanish on $x,y$ and is determined by
its value $ \beta (g)$ which can be any
$p$-th root of unity in $k$ (including $1$, in
which case $ \beta = \varepsilon $ is the
counit of $H$). 
It
follows that 
$$	
	T(h) := \beta (h_{(1)}) l h_{(2)} l^{-1}
\beta ^{-1} (h_{(3)}) 
$$
is the automorphism of $H$ determined by  
$$
	T(g) = g,\qquad
	T(x) = q^i \beta (g)^{-1} x,\qquad
	T(y) = q^{-is} \beta (g)^{-s}y. 
$$
Comparing this with  
$$
	S^2 (g) = g,\qquad
	S^2 (x) = g x g^{-1} = qx,\qquad
	S^2(y) = 
	g^{s} y g^{-s} =
	q^{-s^2} y
$$
shows that $S^2=T$ if and only if 
$$
	\beta (g)  = q^{i-1},\qquad
	\beta (g)^{s} = q^{(s-i)s}.
$$
Assuming $\beta(g)=q^{i-1}$, we obtain 
$$ 
	\beta(l)= \beta (g)^i = q^{i(i-1)}.
$$ 
Furthermore, $\beta (g)^{s} = q^{(s-i)s}$
reduces to $q ^{(1-2i+s)s} = 1$.
In other words, we need that 
$$
	(1-2i+s)s = 0 \in \mathbb{Z} _p.
$$
In total we see that $(l,\beta)$ is a
modular pair in involution if and only if
$$
	(1-2i+s)s = i(i-1) = 0 \in \mathbb{Z} _p.
$$
For $i=0$ this means $s(1+s)=0$ and for $i=1$
it means $s(1-s)=0$ in $ \mathbb{Z} _p$.
The claim follows. \qed


\begin{thebibliography}{99}
\bib{as}{article}{
   author={Andruskiewitsch, Nicol\'{a}s},
   author={Schneider, Hans-J\"{u}rgen},
   title={Hopf algebras of order $p^2$ and braided Hopf algebras of order
   $p$},
   journal={J. Algebra},
   volume={199},
   date={1998},
   number={2},
   pages={430--454},
   issn={0021-8693},
   review={\MR{1489920}},
   doi={10.1006/jabr.1997.7175},
}
\bib{cm}{article}{
   author={Connes, Alain},
   author={Moscovici, Henri},
   title={Cyclic cohomology and Hopf algebras},
   note={Mosh\'{e} Flato (1937--1998)},
   journal={Lett. Math. Phys.},
   volume={48},
   date={1999},
   number={1},
   pages={97--108},
   issn={0377-9017},
   review={\MR{1718047}},
   doi={10.1023/A:1007527510226},
}

\bib{hkrs}{article}{
   author={Hajac, Piotr M.},
   author={Khalkhali, Masoud},
   author={Rangipour, Bahram},
   author={Sommerh\"{a}user, Yorck},
   title={Stable anti-Yetter-Drinfeld modules},
   language={English, with English and French summaries},
   journal={C. R. Math. Acad. Sci. Paris},
   volume={338},
   date={2004},
   number={8},
   pages={587--590},
   issn={1631-073X},
   review={\MR{2056464}},
   doi={10.1016/j.crma.2003.11.037},
}

\bib{kr}{article}{
   author={Kauffman, Louis H.},
   author={Radford, David E.},
   title={A necessary and sufficient condition for a finite-dimensional
   Drinfel\cprime d double to be a ribbon Hopf algebra},
   journal={J. Algebra},
   volume={159},
   date={1993},
   number={1},
   pages={98--114},
   issn={0021-8693},
   review={\MR{1231205}},
   doi={10.1006/jabr.1993.1148},
}

\bib{m}{book}{
   author={Montgomery, Susan},
   title={Hopf algebras and their actions on rings},
   series={CBMS Regional Conference Series in Mathematics},
   volume={82},
   publisher={Published for the Conference Board of the Mathematical
   Sciences, Washington, DC; by the American Mathematical Society,
   Providence, RI},
   date={1993},
   pages={xiv+238},
   isbn={0-8218-0738-2},
   review={\MR{1243637}},
   doi={10.1090/cbms/082},
}

\bib{r}{book}{
   author={Radford, David E.},
   title={Hopf algebras},
   series={Series on Knots and Everything},
   volume={49},
   publisher={World Scientific Publishing Co. Pte. Ltd., Hackensack, NJ},
   date={2012},
   pages={xxii+559},
   isbn={978-981-4335-99-7},
   isbn={981-4335-99-1},
   review={\MR{2894855}},
}
\end{thebibliography}
\end{document}